\documentclass[12pt]{article}
\usepackage{amsmath,amssymb}
\baselineskip 5.5pt \lineskip 5.5pt \numberwithin{equation}{section}

\def \C{\hbox{$C\hskip -5pt \vrule height 6pt depth 0pt \hskip 6pt$}}

\def\qed{\ \ \ifhmode\unskip\nobreak\fi\ifmmode\ifinner
         \else\hskip5pt\fi\fi
 \hbox{\hskip5pt\vrule width4pt height6pt depth1.5pt\hskip 1 pt}}

\def\Span{{\rm span}}

\def\cl{\centerline}

\def\bs{\backslash}

\def\vs{\vspace*}

\def\C{\mathbb{C}}
\def\F{\mathbb{F}}

\begin{document}
\cl {{\Large\bf Verma modules over the generalized}}
 \cl {{\Large\bf
 \vs{10pt} Heisenberg-Virasoro algebra}\footnote{Supported by NSF grants 10471096, 10571120 of China
and ``One Hundred Talents Program'' from University of Science and
Technology of China}} \cl{Ran Shen$^*$,  \ Yucai
Su$^\dag$}\vs{4pt}\cl{\small\it $^*$Department of Mathematics,
Shanghai Jiao Tong University}\cl{\small\it Shanghai 200240,
China}\vs{4pt}\cl{\small\it \cl{\small\it $^\dag$Department of
Mathematics,
 University of Science and Technology of China}} \cl{\small\it Hefei 230026, China}
\vs{4pt}\cl{\small\it Email: ranshen@sjtu.edu.cn, \
ycsu@ustc.edu.cn}\vs{10pt}

{\small {\bf Abstract. } For any additive subgroup $G$ of  an
arbitrary field $\F$ of characteristic zero, there corresponds a
generalized Heisenberg-Virasoro algebra ${\cal L}[G]$. Given a total
order of $G$ compatible with its group structure, and any
$h,h_I,c,c_I,c_{LI}\in \F$, a Verma module
$\widetilde{M}(h,h_I,c,c_I,c_{LI})$ over ${\cal L}[G]$ is defined.
In the this note, the irreducibility of Verma modules $
\widetilde{M}(h,h_I,c,c_I,c_{LI})$ is completely determined.
 \vs{2pt}\par {\bf Key
Words:} The generalized Heisenberg-Virasoro algebra, Verma modules
\vs{2pt}\par{\it  Mathematics Subject Classification (2000)}: 17B56;
17B68.}

\vs{10pt}\par \cl{\bf1. \ Introduction}
\par
\def\C{{\mathbb{F}}}
Let $\F$ be a field of characteristic 0. The well-known {\it twisted
Heisenberg-Virasoro algebra} is the Lie algebra ${\cal L}:={\cal
L}[\mathbb{Z}]$ with an $\F$-basis $\{L_m,I_m,C,C_I,C_{LI}\,|\, m
\in \mathbb{Z}\}$ subject to the following relations (e.g., [ACKP,
B]) \begin{eqnarray*} \!\!\!&\!\!\!&
[L_n,L_m]=(n-m)L_{n+m}+\delta_{n,-m}\frac{n^{3}-n}{12}C,
\\[4pt]\!\!\!&\!\!\!&
[L_n,I_m]=-mI_{n+m}-\delta_{n,-m}(n^{2}+n)C_{LI}, \\[4pt]\!\!\!&\!\!\!&
[I_n,I_m]=n\delta_{n,-m}C_{I},\\[4pt] \!\!\!&\!\!\!& [{\mathcal
L},C]=[{\mathcal L},C_{LI}]=[{\mathcal L},C_{I}]=0.
\end{eqnarray*}
This Lie algebra is the universal central extension of the Lie
algebra of differential operators on a circle of order at most one,
which contains an infinite-dimensional Heisenberg subalgebra and the
Virasoro subalgebra. The natural action of the Virasoro subalgebra
on the Heisenberg subalgebra is twisted with a 2-cocycle. The
structure and representation theory for the twisted
Heisenberg-Virasoro algebra has been well developed (e.g., [ACKP, B,
FO, JJ,  SJ]). The structure of the irreducible highest weight
modules for the twisted Heisenberg-Virasoro algebra are determined
in [ACKP, B].

By replacing the index group $\mathbb{Z}$ by an arbitrary subgroup
$G$ of the base field $\F$, it is natural to introduce the so-called
{\it generalized Heisenberg-Virasoro algebra ${\cal L}[G]$}
(cf.~Definition 2.1, see e.g., [XLT, LJ]). This is the Lie algebra
which is the 3-dimensional universal central extension of the Lie
algebra of generalized differential operators of order at least one.
The Harish-Chandra modules of intermediate series over generalized
Heisenberg-Virasoro algebra ${\cal L}[G]$ are discussed in [LJ].

Given any total order of $G$ compatible with its group structure,
and given any $h,h_I,c,c_I,$ $c_{LI}\in \F$, there corresponds a
Verma module $\widetilde{M}(h,h_I,c,c_I,c_{LI})$ over ${\cal L}[G]$.
Due to the fact that the representations of generalized
Heisenberg-Virasoro algebras are closely related to the
representation theory of toroidal Lie algebras as well as some
problems in mathematical physics (e.g., [ACKP, FO, JJ]) and the
Verma modules play the crucial role in the representation theory, it
is very natural to consider the Verma modules over the generalized
Heisenberg-Virasoro algebras. In this note, we completely determine
the irreducibility of Verma modules
$\widetilde{M}(h,h_I,c,c_I,c_{LI})$ over ${\cal L}[G]$ for arbitrary
$G$. Namely, if $G$ does not contain a minimal positive element with
respect to the total order, then the Verma module
$\widetilde{M}(h,h_I,c,c_I,c_{LI})$ is irreducible if and only if
$(c_I,c_{LI})\neq (0,0)$; in case if $G$ contains the minimal
positive element $a$, then the  Verma module
$\widetilde{M}(h,h_I,c,c_I,c_{LI})$ is irreducible if and only if
the ${\cal L}[\mathbb{Z}a]$-module generated by a fixed highest
weight generator is irreducible over the twisted Heisenberg-Virasoro
algebra ${\cal L}[\mathbb{Z}a]$ (cf.~Theorem 3.1).
 \vs{10pt}
\par

\cl{\bf2. \ Generalized Heisenberg-Virasoro algebras}
\setcounter{section}{2}\setcounter{theo}{0} \setcounter{equation}{0}
\par
Let $U:=U({\cal L})$ be the universal enveloping algebra of the
twisted Heisenberg-Virasoro algebra ${\cal L}$. For any
$h,h_I,c,c_I,c_{LI}\in \F$, denote by $I(h,h_I,c,c_I,c_{LI})$ the
left ideal of $U$ generated by the elements
$$\{L_i,I_j\,|\,i,j>0\}\cup\{L_0-h\cdot 1,I_0-h_I\cdot 1,C-c\cdot 1,C_I-c_I\cdot 1,C_{LI}-c_{LI}\cdot
1\}.$$ The {\it Verma module} with highest weight
$(h,h_I,c,c_I,c_{LI})$ over $\cal L$ is defined as
$$
M(h,h_I,c,c_I,c_{LI}):=U/I(h,h_I,c,c_I,c_{LI}),$$ which is a highest
weight module  with a basis consisting of all vectors of the form
\begin{equation}I_{-p_1}I_{-p_{2}}\cdots I_{-p_s}L_{-j_1}L_{-j_2}\cdots
L_{-j_k}v_h,\end{equation} where $ s,k\in \mathbb{N}\cup\{0\},\,
p_r,j_i\in \mathbb{N}$ and $0< p_1\leq p_{2}\leq\cdots \leq p_s,\,0<
j_1\leq j_2\leq\cdots\leq j_k.$ \vs{4pt}\par {\bf Definition 2.1 }
Let $G\subseteq \F$ be an additive subgroup. The {\it generalized
Heisenberg-Virasoro algebra} $\widetilde{\cal L}:={\cal L}[G]$ is a
Lie algebra with $\F$-basis $\{L_\mu,I_\mu,C,C_I,C_{LI}\,|\, \mu\in
G\}$ subject to the following relations [XLT, LJ]
\begin{eqnarray*}
\!\!\!&\!\!\!&
[L_{\mu},L_{\nu}]=(\mu-\nu)L_{\mu+\nu}+\delta_{\mu,-\nu}\frac{\mu^{3}-\mu}{12}C,\\[4pt]
\!\!\!&\!\!\!& [L_{\mu},I_{\nu}]=-\nu
I_{\mu+\nu}-\delta_{\mu,-\nu}(\mu^{2}+\mu)C_{LI},\\[4pt]
\!\!\!&\!\!\!& [I_\mu,I_\nu]=\mu\delta_{\mu,-\nu}C_{I},\\[4pt]
\!\!\!&\!\!\!& [\widetilde{\mathcal L},C]=[\widetilde{\mathcal
L},C_{LI}]=[\widetilde{\mathcal L},C_{I}]=0.
\end{eqnarray*}

For any $x\in G^*:=G\setminus\{0\}$, obviously,
$\mathbb{Z}x\subseteq G$. Let ${\cal L}[\mathbb{Z}x]$ be the
$\F$-subspace of $\widetilde{\mathcal L}$ spanned by
$\{L_{ix},I_{ix},C,C_I,C_{LI}\,|\, i\in \mathbb{Z}\}$. It is clear
that
 ${\cal L}[\mathbb{Z}x]$ is a Lie algebra
isomorphic to the twisted Heisenberg-Virasoro algebra $\cal L$.
Precisely, we have \vs{4pt}\par {\bf Lemma 2.2 }  The map
\begin{equation*}\begin{array}{lllll}\theta:&{\cal L} &
\rightarrow & {\cal L}[\mathbb{Z}x]\\[4pt] &L_i & \mapsto &
x^{-1}L_{ix}+\delta_{i,0}\frac{x-x^{-1}}{24}C,\\[4pt] &I_i&\mapsto&
x^{-1}I_{ix}+\delta_{i,0}(1-x^{-1})C_{LI},\\[4pt]
&C & \mapsto & xC,\\[4pt] &C_I&\mapsto& x^{-1}C_I,\\[4pt] &C_{LI}&\mapsto&
C_{LI},
\end{array}\end{equation*}
for $i\in \mathbb{Z},$ extends uniquely to a Lie algebra isomorphism
between ${\cal L}$ and ${\cal L}[\mathbb{Z}x]$. \vs{4pt}\par{\it
Proof.~}~This follows from straightforward
verifications.$\hfill\Box$ \vs{4pt}\par Throughout this note, we fix
a total order ``$\succ$'' on $G$ compatible with its group
structure, namely, $x\succ y$ implies $x+z\succ y+z$ for any $z\in
G$. Denote
$$G_+:=\{x\in G\,|\,x\succ 0\},\ \ \ G_-:=\{x\in G\,|\,x\prec 0\}.$$ Then
$G=G_+\cup\{0\}\cup G_-$.

For an $\widetilde{\mathcal L}$-module $V$ and
$\lambda,h_I,c,c_I,c_{LI}\in \F$, denote by
$$V_{\lambda,h_I,c,c_I,c_{LI}}:=\{v\in V|L_0v=\lambda v,\ I_0v=h_Iv,\
Cv=cv,\ C_Iv=c_Iv,\ C_{LI}v=c_{LI}v\},$$ the {\it weight space} of
$V$. We shall simply write $V_\lambda$ instead of
$V_{\lambda,h_I,c,c_I,c_{LI}}$. Define
$$\mbox{supp$(V):=\{\lambda\in \F\,|\,V_\lambda\neq 0\}$,}$$
called the {\it weight set} (or the {\it support}) of $V$. For any
$h,h_I,c,c_I,c_{LI}\in \F$, let $\widetilde{M}(h,h_I,c,c_I,c_{LI})$
be the {\it Verma module} for $\widetilde{\mathcal L}$, which is
defined by using the order ``$\succ$'' and the same fashion as that
for $\cal L$ at the beginning of this section. Then $I_0,\ C,\ C_I,\
C_{LI}$ acts as $h_I,c,c_I,c_{LI}$ respectively on
$\widetilde{M}(h,h_I,c,c_I,c_{LI})$ and
$$\mbox{supp$(\widetilde{M}(h,h_I,c,c_I,c_{LI}))=h+G_+$.}$$
For any $x\in G_+$, let
$$\widetilde{M}_x(h,h_I,c,c_I,c_{LI})=
U({\cal L}[\mathbb{Z}x])v_h,\eqno(2.2)$$ be the ${\cal
L}[\mathbb{Z}x]$-submodule of $\widetilde{M}(h,h_I,c,c_I,c_{LI})$
generated by a fixed highest weight generator $v_h$. Note that the
subgroup $\mathbb{Z}x$ is also a ``totally ordered abelian group'',
inheriting the order ``$\succ$'' from $G$. It is easy to see that
$$ax\succ bx \ \ \ \Longleftrightarrow\ \ \ a>b \ \mbox{ for }\ a,b\in \mathbb{Z}.$$
As a result, we have \vs{4pt}\par {\bf Corollary 2.3~}~As an ${\cal
L}$-module, we have
$$\widetilde{M}_x(h,h_I,c,c_I,c_{LI})\cong
M(x^{-1}h+\frac{x-x^{-1}}{24}c,x^{-1}h_I+(1-x^{-1})c_{LI},xc,x^{-1}c_I,c_{LI}).$$
\vs{4pt}\par {\it Proof.~}~This is clear by Lemma 2.2.$\hfill\Box$
\vskip10pt

\cl{\bf3. \ The main result}
\setcounter{section}{3}\setcounter{theo}{0} \setcounter{equation}{3}
\par

Recall that $(G,\succ)$ is a totally ordered abelian group. Denote
$$\mbox{$B(x)=\{y\in G\,|\,0\prec y\prec x\}$ \ \ for $x\in G_+$.}$$ The order
``$\succ$'' is called {\it dense} if $\sharp B(x)=\infty$ for all
$x\in G_+$; {\it discrete} if there exists some $a\in G_+$ such that
$B(a)=\emptyset$, in this case $a$ is called the {\it minimal
positive element} of $G$.

For convenience, we denote \vs{6pt} \\ \
\hspace*{20pt}$L_{-j}:=L_{-j_1}L_{-j_2}\cdots L_{-j_k}\mbox{ \ \ for
\ } 0\prec j_1\preceq j_2\preceq\cdots\preceq j_k,\
j=(j_1,j_2,\cdots,j_k),$\hfill(3.1) \vs{6pt}\\ \ \hspace*{20pt}$
I_{-p}:=I_{-p_s}I_{-p_{s-1}}\cdots I_{-p_1}\mbox{ \ \ for \ }0\prec
p_s\preceq \cdots\preceq p_2\preceq p_1,\
p=(p_s,\cdots,p_2,p_1).$\hfill(3.2) \vs{6pt}\\
Then $U(\widetilde{\mathcal L}_-)$ has a basis
$$\{I_{-p}L_{-j}\ |\mbox{ \ for all \ $j,p$ as in (3.1) and (3.2)}\}. \eqno(3.3)$$
 Denote by $|j|$ the number of components in $j$. Then
$|j|=k$ in (3.1) and $|p|=s$ in (3.2).

The main result in this note is following. \vs{4pt}\par {\bf Theorem
3.1~}~{\it Let $h,h_I,c,c_I,c_{LI}\in \F$\vs{-4pt}.
\begin{itemize}\parskip-6pt
\item[{\rm(1)}] With respect to a dense order ``$\succeq$'' of $G$,
the Verma module $\widetilde{M}(h,h_I,c,c_I,c_{LI})$ is an
irreducible ${\cal L}[G]$-module if and only if $(c_I,c_{LI})\neq
(0,0)$.
\item[{\rm(2)}] With respect to a discrete order
``$\succeq$'' of $G$ with minimal positive element $a$, the Verma
module $\widetilde{M}(h,h_I,c,c_I,c_{LI})$ is an irreducible ${\cal
L}[G]$-module if and only if $\widetilde{M}_a(h,h_I,c,c_I,c_{LI})$
$($cf.~$(2.2))$ is an irreducible ${\cal L}[\mathbb{Z}a]$-module.
\end{itemize}
}\par {\bf Remark 3.2~}~{\it  Suppose $c_I=c_{LI}=0$ in case of
Theorem $3.1(1)$. Since
$$\widetilde{\mathcal I}:=\Span_\F\{I_\mu,C_I,C_{LI}\,|\,\mu\in G\},$$ is an ideal
of $\widetilde{\mathcal L}$, the Verma module
$V:=\widetilde{M}(h,h_I,c,0,0)$ over $\widetilde{\mathcal L}$ has a
proper submodule $U(\widetilde{\mathcal I})V$ such that the quotient
module $W:=V/U(\widetilde{\mathcal I})V$ is simply the Verma module
over the generalized Virasoro algebra ${\rm
Vir}[G]:=\Span_\F\{L_\mu,C\,|\,\mu\in G\}\cong \widetilde{\mathcal
L}/\widetilde{\mathcal I}$, whose irreducibility is completely
determined in {\rm[HWZ].} Also note that the irreducibility of a
Verma module over the twisted Heisenberg-Virasoro algebra $\mathcal
L$ is completely determined in {\rm [B]}. Thus, essentially the
above theorem has in fact determined the structure of all Verma
modules over $\widetilde{\mathcal L}$.}\vs{4pt}

{\it Proof of Theorem 3.1.~}~(1) Suppose the order ``$\succeq$'' of
$G$ is dense. Let $v_h$ be a fixed highest weight generator in
$\widetilde{M}(h,h_I,c,c_I,c_{LI})$ of weight $h$. Let $u_0\not\in
\F v_h$ be any given weight vector in
$V:=\widetilde{M}(h,h_I,c,c_I,c_{LI})$. \vskip4pt\par {\bf Claim 1}:
There exists a weight vector $u\in U({\cal L}[G])u_0$ of weight
$\lambda$ such that
$$u=\sum\limits_p a_pI_{-p}v_h\mbox{ (a finite sum) \
for some }a_p\in\F^*=\F\bs\{0\}.\eqno(3.4)$$

For each $m \in \mathbb{N}$, set
$$V_m:=\sum\limits_{p,j:\,|j|\leq m}\F I_{-p}L_{-j}v_h.\eqno(3.5)$$
It is clear that $$\mbox{$L_xV_m\subseteq V_m,\ I_xV_m\subseteq V_m$
\ for  $x\in G_+$.}$$ We can write $u_0$ as (cf.~(2.1) and (3.3))
$$u_0=\sum\limits_{p,j}a_{pj}I_{-p}L_{-j}v_h\mbox{ \ for some }a_{pj}\in\F^*.$$
 Let $r:=$ max $\{|j|\ |\,a_{pj}\neq
0\}$. If $r=0$, then the claim holds clearly. We assume $r\geq 1$,
and write
$$u_0\equiv u_0' \,(\mbox {mod\,} V_{r-1}),\
\mbox{ where }u_0'=\sum\limits_{p,j:\,|j|=r}a_{pj}I_{-p}L_{-j}v_h.
\eqno(3.6)$$ Let $x\in G_+$ such that (cf.~(3.2) for notation $p_l$)
$$\mbox{ $x\prec
\mbox{min}\{j_1\,|\,a_{pj}\neq 0\}$ and $\{x,j_1-x\,|\,a_{pj}\neq
0\}\cap \{p_l\,|\,a_{pj}\neq 0,\ \forall\, l\}=\emptyset.$}$$ Then
$$I_xu_0'=\sum\limits_{p,j:\,|j|=r}xa_{pj}I_{-p}\left(\sum\limits_{i=1}^rL_{-j_1}
\cdots L_{-j_{i-1}}I_{x-j_i}L_{-j_{i+1}}\cdots L_{-j_r}\right)v_h.$$
If any $$\mbox{$xa_{pj}I_{-p}L_{-j_1} \cdots
L_{-j_{i-1}}I_{x-j_i}L_{-j_{i+1}}\cdots L_{-j_r}$ and
$xa_{p'j'}I_{-p'}L_{-j_1'} \cdots
L_{-j_{s-1}'}I_{x-j_s'}L_{-j_{s+1}'}\cdots L_{-j_r'}$,}$$ for $1\leq
i,s\leq r,$  are linear dependent, it is not difficult to see that
$p=p'$ and $j=j'$. Hence $$0\neq u_1:=I_xu_0'\in V_{r-1}.$$
Similarly, let $u_1 \equiv u_1'\ (\mbox{mod }V_{r-2})$ as in (3.6),
then $u_1'\neq 0$. For $k=2,\cdots,r$. We define recursively and
prove by induction that,
$$u_k: = I_xu_{k-1}\in V_{r-k},\ u_k\equiv u_k'\ (\mbox{mod
}V_{r-k-1}), \ u_k'\neq 0.$$Letting $k=r$, we get that $0\neq u_r\in
V_0$. Our claim follows. \vskip4pt Now let $u$ be as in (3.4). Set
$P:=\{p\,|\,a_p\neq 0\}\neq\emptyset$. We define the total order
``$\succ$'' on $P$ as follows: For any $p,p'\in P$, if
$k:=|p|>l:=|p'|$, we set $p'_i=0$ for $i=l+1,...,k$. Then
\begin{equation*}
p\succ p'\ \ \ \Longleftrightarrow\ \ \ \exists\,s\mbox{ with }
1\leq s\leq k\text{ such that }p_s\succ p_s'\text{ and }p_t=p_t'
\text{ for  }t<s. \eqno(3.7)\end{equation*} Let
$$q:=(q_{k_0},\cdots,q_2,q_1),\ \ 0\prec q_{k_0}\preceq\cdots\preceq
q_1,$$ be the unique maximal element in $P$. Then \vskip4pt Case 1:
If $c_I\neq 0,\ $ then by the simple calculations
$$bv_h=I_qu\in U({\cal L}[G])u_0\mbox{ \ for some \ }b\in \F ^*.$$

Case 2: Suppose $c_I=0,\ c_{LI}\neq 0$. Let $y\in G_+$ such that
$$\{x\in G\,|\,q_1-y\prec x\prec q_2\}\cap\{p_1,p_2\,|\,p\in
P\}=\emptyset.$$
 Then
$$u':=L_{q_1-y}u=a'I_{-z}v_h
\mbox{ \ for some \ }a'\in\F^*,$$ where
$$\begin{array}{ll}
z=(z_{k_0},\cdots,z_2,z_1),\ \ 0\prec z_{k_0}\preceq\cdots\preceq
z_2\preceq z_1,\mbox{ \ and}\\[6pt]
\{z_i\,|\,i=1,2,\cdots,k_0\}=\{q_{k_0},\cdots,q_3,q_2,y\}.\end{array}$$

(i) If $\{z_i\,|\,i=1,2,\cdots,k_0\}\cap\{h_I/c_{LI}-1\}=\emptyset,$
then
$$b'v_h=L_zu'\in U({\cal L}[G])u_0\ne0,\mbox{ where } b'=\prod\limits_{i=1}^{k_0}z_i(h_I-(z_i+1)c_{LI})\in \F ^*.$$

(ii) If there exists some $z_i=h_I/c_{LI}-1$ with $1\leq i\leq k_0.$
We assume $$\{z_i\}\cap\{z_k\,|\,1\leq k\leq k_0, k\neq
i\}=\emptyset.$$ Otherwise, we only need to recurse the following
proof. Let
$$w:=L_{z_{i-1}}\cdots L_{z_2}L_{z_1}u'=a''I_{-z_{k_0}}I_{-z_{k_0-1}}\cdots I_{-z_i}v_h\neq 0
\mbox{ for some }a''\in\F^*.$$ Take $x'\in G_+$ such that
$z_i-x'\succ z_k, i<k\leq k_0$. Then
$$w':=L_{z_i-x'}w=\overline{a}I_{-z_{k_0}}I_{-z_{k_0-1}}\cdots I_{-z_{i+1}}I_{-x'}v_h\neq 0\mbox{ for
some }\overline{a}\in\F^*,$$ and
$\{z_{k_0},z_{k_0-1}\cdots,z_{i+1},x'\}\cap\{h_I/c_{LI}-1\}=\emptyset$.
This becomes case (i) if we take $u'$ to be $w'$.

Therefore, $v_h\in U({\cal L}[G])u_0$ in any case. Hence
$\widetilde{M}(h,h_I,c,c_I,c_{LI})$ is irreducible. \vskip4pt (2)
Suppose the order ``$\succ$'' of $G$ is discrete with the minimal
positive element $a$. Then $\mathbb{Z} a\subseteq G$. For any $x\in
G$, we write $x\succ \mathbb{Z}a$ if $x\succ na$ for all $n\in
\mathbb{Z}$. Let $$H_+:=\{x\in G\,|\,x\succ \mathbb{Z}a\},\ \
H_-=-H_+.$$ It is not difficult to see that $$G=\mathbb{Z}a\cup
H_+\cup H_-.\eqno(3.8)$$ Then one can see that $$\mbox{${\cal
L}[H_+]\widetilde{M}_a(h,h_I,c,c_I,c_{LI})=0$ \ (recall (2.2)).}$$
Since
$$\widetilde{M}(h,h_I,c,c_I,c_{LI})\cong U({\cal
L}[G])\otimes_{U({\cal L}[\mathbb{Z}a]+{\cal
L}[H_+])}\widetilde{M}_a(h,h_I,c,c_I,c_{LI}),$$ it follows that the
irreducibility of ${\cal L}[G]$-module
$\widetilde{M}(h,h_I,c,c_I,c_{LI})$  imply the irreducibility of
${\cal L}[\mathbb{Z}a]$-module
$\widetilde{M}_a(h,h_I,c,c_I,c_{LI})$.

Conversely, suppose $\widetilde{M}_a(h,h_I,c,c_I,c_{LI})$ is an
irreducible ${\cal L}[\mathbb{Z}a]$-module. Let $u_0\not\in \F v_h$
be any weight vector in $\widetilde{M}(h,h_I,c,c_I,c_{LI})$. We want
to prove $$U({\cal L}[G])u_0\cap
\widetilde{M}_a(h,h_I,c,c_I,c_{LI})\neq \{0\},\eqno(3.9)$$ from
which the irreducibility of $\widetilde{M}(h,h_I,c,c_I,c_{LI})$ as
${\cal L}[G]$-module follows immediately. \vskip4pt Case 1: $c_I\neq
0.$ We can write $u_0$ as (cf.~(3.8))
$$u_0\equiv\sum\limits_{p_l',j_k'\in H_+,\,p_s,j_r\in
\mathbb{Z}_+a,\,|j'|+|j|=r}a_{p'j'pj}I_{-p'}L_{-j'}I_{-p}L_{-j}v_h\
(\mbox{mod }V_{r-1})\mbox{ \ for some }a_{p'j'pj}\in\F^*,$$ where
$V_{r-1},\,r$ are defined as in (3.5) and (3.6).
 Let (cf.~notation (3.2))
$$
\begin{array}{lll}
P'=\{pp'=(p_s,\cdots,p_2,p_1,p_l',\cdots,p_2',p_1')\,|\!\!\!\!&0\prec
p_s\preceq \cdots\preceq p_2\preceq p_1\prec \\[4pt]&
p_l'\preceq \cdots\preceq p_2'\preceq p_1',\ \, a_{p'j'pj}\neq
0\}.\end{array}$$ If $P'\ne\emptyset$, we define the total order
``$\succ$'' on $P'$ as in (3.7). Let $q^0$ be the maximal element in
$P'$. Then
$$u_0':=I_{q^0}u_0\equiv\sum\limits_{j_k'\in
H_+,\,j_r\in \mathbb{Z}_+a,\,|j'|+|j|=r}a_{j'j}L_{-j'}L_{-j}v_h\
(\mbox{mod }V_{r-1})\mbox{ \ for some }a_{j'j}\in\F^*.$$ If
$P'=\emptyset$, then $u_0$ has the form of $u_0'$ naturally. By the
proof of [HWZ, Theorem 3.1], there exists a weight vector $0\ne u\in
U({\cal L}[G])u_0\cap \widetilde{M}_a(h,h_I,c,c_I,c_{LI})$, which
gives (3.9) as required.

Case 2: $c_I=0.$ We can write$$u_0=\sum\limits_{p_l',j_k'\in
H_+,\,p_s,j_r\in
\mathbb{Z}_+a}b_{p'j'pj}I_{-p'}L_{-j'}I_{-p}L_{-j}v_h\mbox{ \ for
some }b_{p'j'pj}\in\F^*.$$
 If $J:=\{j'\,|\,b_{p'j'pj}\neq
0\}\neq\emptyset$, we set $j(0):=\mbox{min}\{j_1'\,|\,b_{p'j'pj}\neq
0\}$. Then there exists some $m\in \mathbb{N}$ such that
$$\{j_1'-\varepsilon\,|\,b_{p'j'pj}\neq
0\}\cap\{p_l'\,|\,b_{p'j'pj}\neq 0, \,\forall\, l\}=\emptyset,
\mbox{ \ where }\varepsilon=j(0)-ma.$$ Let $n_0=\mbox{max }\{|j'|\
|\,b_{p'j'pj}\neq 0\},$ then
$$u':=I_\varepsilon^{n_0}u_0=\sum\limits_{p_l'\in
H_+,\,p_s,j_r\in \mathbb{Z}_+a}b_{p'pj}'I_{-p'}I_{-p}L_{-j}v_h\neq
0\mbox{ \ for some }b_{p'pj}'\in\F^*,$$by the proof of Claim 1. If
$J=\emptyset$, then $u_0$ has the form of $u'$ naturally. Let
$$Q:=\{p'\,|\,b_{p'pj}'\neq 0,\
|p'|=t\},\mbox{ \ where }t=\mbox{min}\{|p'|\ |\,b_{p'pj}'\neq 0\}.$$
If $t=0$, the theorem holds clearly since $u'$ is a weight vector.
We assume $t\geq 1.$ Then $Q\neq\emptyset$. Again, we define the
total order ``$\prec$'' on $Q$ as in (3.7). Let
$$q':=(q_1',q_2',\cdots,q_t'),\ \ 0\prec q_1'\preceq q_2'\preceq
\cdots\preceq q_t',$$ be the unique minimum element in $Q$. For
$m\in \mathbb{N}$, set$$V_m':=\sum\limits_{p_l'\in H_+,\,p_s,j_r\in
\mathbb{Z}_+a,\,|p'|\geq m}\F I_{-p'}I_{-p}L_{-j}v_h.$$Then
$$u'\equiv\sum\limits_{p_l'\in H_+,p_s,j_r\in
\mathbb{Z}_+a,|p'|=t}b_{p'pj}'I_{-p'}I_{-p}L_{-j}v_h\
(\mbox{mod}V_{t+1}').$$We
have$$u(1):=L_{q_1'-a}u'\equiv\sum\limits_{p_l^{(1)}\in
H_+,\,p_s,j_r\in
\mathbb{Z}_+a,\,|p^{(1)}|=t-1}b_{p^{(1)}pj}^{(1)}I_{-p^{(1)}}I_{-a}I_{-p}L_{-j}v_h\
(\mbox{mod}V_t')$$ for some $b_{p^{(1)}pj}^{(1)}\in\F^*$. Define
$Q^{(1)}=\{p^{(1)}\,|\,b_{p^{(1)}pj}^{(1)}\neq 0\},\
q^{(1)}=(q_2',q_3',\cdots,q_t').$ By our assumption and the
commutator relations for ${\cal L}[G]$, we see that
$b_{q^{(1)}pj}^{(1)}\neq 0$, hence $Q^{(1)}\neq\emptyset$. Moreover,
$q^{(1)}$ is the unique minimum element in $Q^{(1)}$.

Now for $s=2,3,\cdots,t$, we define recursively and prove by
induction that

(i) Let $u(s):=L_{q_s'-a}u(s-1)$. Then
$$u(s)\equiv\sum\limits_{p_l^{(s)}\in H_+,\,p_s,j_r\in
\mathbb{Z}_+a,\,|p^{(s)}|=t-s}b_{p^{(s)}pj}^{(s)}I_{-p^{(s)}}I_{-a}^sI_{-p}L_{-j}v_h\
(\mbox{mod}V_{t-s+1}')\mbox{ \ for some
}b_{p^{(s)}pj}^{(s)}\in\F^*.$$

(ii) Let $Q^{(s)}:=\{p^{(s)}\,|\,b_{p^{(s)}pj}^{(s)}\neq
0\}\neq\emptyset$. Moreover,
$q^{(s)}:=(q_{s+1}',q_{s+2}',\cdots,q_t')$ is the unique minimum
element in $Q^{(s)}$.

Now letting $s=t$ and noting that $u(t)$ is a weight vector, we get
that $0\neq u(t)\in U({\cal L}[G])u_0\cap
\widetilde{M}_a(h,h_I,c,c_I,c_{LI})$, which gives (3.9) as
required.$\hfill\Box$

\par
 \vs{10pt}
\par
\cl{\bf REFERENCES} \vs{-1pt}\small \lineskip=4pt
\begin{itemize}\parskip-5pt
\item[{[ACKP]}] E. Arbarello, C. De Concini, V.G. Kac, C.
Procesi, Moduli spaces of curves and representation theory, {\it
Comm. Math. Phys.}, {\bf117}(1988), 1-36.

\item[{[B]}] Y. Billig, Respresentations of the twisted
Heisenberg-Virasoro algebra at level zero, {\it Canad. Math.
Bulletin}, {\bf46}(2003), 529-537.

\item[{[FO]}] M.A. Fabbri, F. Okoh, Representations of
Virasoro-Heisenberg algebras and Virasoro-toroidal algebras, {\it
Canad. J. Math.}, {\bf51}(1999), no.3, 523-545.

\item[{[HWZ]}] J. Hu, X. Wang, K. Zhao, Verma modules over
generalized Virasoro algebras $Vir[G]$, {\it J. Pure Appl. Algebra,}
{\bf177}(2003), no.1, 61-69.

\item[{[JJ]}] Q. Jiang, C. Jiang, Representations of the twisted
Heisenberg-Virasoro algebra and the full toroidal Lie algebras, {\it
Algebra Colloq.}, accepted.


\item[{[LJ]}] D. Liu, C. Jiang, The generalized Heisenberg-Virasoro algebra, preprint
(arXiv:math.RT/ 0510545).

\item[{[SJ]}] R. Shen, C. Jiang, Derivation algebra and automorphism
group of the twisted Heisenberg-Virasoro algebra, preprint.

\item[{[WZ]}] X. Wang, K. Zhao, Verma modules over
the Virasoro-like algebra, {\it J. Aust. Math. Soc.,} in press.

\item[{[XLT]}]
M. Xue, W. Lin, S. Tan, Central extension, derivations and
automorphism group for Lie algebras arising from the 2-dimensional
torus, {\it Journal of Lie Theory}, {\bf16}(2005), 139-153.

\end{itemize}
\end{document}